\documentclass[reqno]{amsart}
\newtheorem{Pa}{Paper}[section]

\newtheorem{theorem}[Pa]{{\bf Theorem}}
\newtheorem{Tm}[Pa]{{\bf Theorem}}
\newtheorem{La}[Pa]{{\bf Lemma}}

\newtheorem{Rk}[Pa]{{\bf Remark}}

\newtheorem{Ex}[Pa]{{\bf Example}}

\newcommand{\C}{{\mathbb C}}
\newcommand{\D}{{\mathbb D}}

\newcommand{\M}{{\mathcal T}}

\newcommand{\cE}{{\mathcal E}}
\newcommand{\cS}{{\mathcal S}}

\newcommand{\cR}{{\mathcal R}}

\newcommand{\cSR}{\mathcal{SR}}  
\newcommand{\bal}{\boldsymbol \alpha}
\newcommand{\bbe}{\boldsymbol \beta}

\begin{document}
\title[Low degree rational solutions]
{An algorithm for finding low degree rational solutions 
to the Schur coefficient problem}

\author{Vladimir Bolotnikov}
\address{Department of Mathematics,
The College of William and Mary,
Williamsburg VA 23187-8795, USA}
\subjclass{41A05, 41A20, 30E05}
\keywords{Schur problem, low degree rational interpolants}

\begin{abstract}
We present an algorithm producing all rational functions 
$f$ with  prescribed $n+1$ Taylor coefficients at the origin
and such that $\|f\|_\infty\le 1$ and $\deg \, f\le k$ for every fixed 
$k\ge n$. The case where $k<n$ is also discussed.
\end{abstract}
\maketitle

\section{Introduction}
\setcounter{equation}{0}

Let $H^\infty$ be the Banach space of bounded analytic functions on the 
open unit disk ${\mathbb D}$ with norm 
$\|f\|_\infty:=\sup_{z\in\D}|f(z)|<\infty$. The closed unit ball $\cS$ of 
$H^\infty$ (sometimes called {\em the Schur class}) thus consists of analytic 
functions mapping $\D$ into its closure. The classical Schur problem which 
we will denote by ${\bf SP}_n$ consists of finding $f\in\cS$  having 
prescribed $n+1$ Taylor coefficients at the origin.

\medskip
${\bf SP}_n$: {\em Given $c_0,\ldots,c_{n}\in\C$, find all functions 
$f\in\cS$ of the form}
\begin{equation}
f(z)=c_0+c_1z+\ldots+c_{n}z^{n}+O(z^{n+1}).
\label{1.1}   
\end{equation}
The problem has a solution if and only if the Pick matrix of the problem
given by 
$$
P_n=I-\M(c_n,\ldots,c_0)\M(c_n,\ldots,c_0)^*
$$
is positive semidefinite. Here and in what follows, $I$ denotes the 
identity matrix of the size  always clear from the context,
and $\M(c_0,\ldots,c_n)$ stands for the lower triangular Toeplitz matrix 
with the bottom row entries indicated in the parentheses:
\begin{equation}
\M(c_n,\ldots,c_0):=
\left[\begin{array}{cccccc} c_0&0&0&\cdots&0\\ c_1
&c_0&0&\cdots&0\\ \vdots&\vdots&\vdots& &0\\ c_{n}&c_{n-1}&
c_{n-2}&\cdots&c_0\end{array}\right].
\label{1.2}
\end{equation}
If $P_n\ge 0$ is singular, then the problem ${\bf SP}_n$ has
a unique solution which is a finite Blaschke product of degree equal to
the rank of $P_n$. In what follows, we assume that the  data set 
$\{c_0,\ldots,c_n\}$ is such that $P_n>0$ and we will call such a data set 
{\em admissible}. For an admissible data set, the parametrization of the 
solution set of the problem ${\bf SP}_n$ was established in 
\cite{Schur} via  the famous Schur algorithm which we now recall.  
Starting with $c_0,\ldots,c_{n}$, define the numbers $c_k^{(j)}$ 
($j=1,\ldots,n; 
k=0,\ldots,n-j$) from  the following recursion:
\begin{equation} 
\left[\begin{array}{c}c_0^{(0)} \\ c_1^{(0)} \\ \vdots \\
c_{n}^{(0)}\end{array}\right]=\left[\begin{array}{c}c_0 \\ c_1 
\\ \vdots \\ c_{n}\end{array}\right]\quad\mbox{and}\quad
\left[\begin{array}{c}c_0^{(j+1)} \\ c_1^{(j+1)} \\ \vdots \\
c_{n-j-1}^{(j+1)}\end{array}\right]=M_j^{-1}
\left[\begin{array}{c}c_1^{(j)} \\ c_2^{(j)} \\ \vdots \\
c_{n-j}^{(j)}\end{array}\right]\quad (j\ge 0),
\label{1.3}
\end{equation}
where the matrix
$$
M_j=\M\left(-\overline{c}_0^{(j)}c_{n-j-1}, \;  \ldots, \; 
-\overline{c}^{(j)}_0c^{(j)}_2, \; -\overline{c}^{(j)}_0c^{(j)}_1, 1-|c^{(j)}_0|^2
\right)
$$
is defined via formula \eqref{1.2}. Let
\begin{equation}
\gamma_j=c_0^{(j)}\quad\mbox{for}\quad j=0,\ldots,n.
\label{1.4}
\end{equation}
If $c_0,\ldots,c_n$ are the Taylor coefficients of an $f\in\cS$, then the 
numbers $\gamma_i$ constructed above are the $n+1$ first Schur parameters 
of $f$ and condition $P_n>0$ is equivalent to $|\gamma_i|<1$ for 
$i=0,\ldots,n$. The Schur algorithm relies on the 
following fact:

\medskip

{\em A function $f$ belongs to $\cS$ and satisfies 
\eqref{1.1} if and only if it is of the form 
\begin{equation}
f(z)=\frac{zf_1(z)+c_0}{z\bar{c}_0f_1(z)+1}
\label{1.5}
\end{equation}
for some $f_1\in\cS$ such that
$\; f_1(z)=c^{(1)}_0+c^{(1)}_1z+\ldots+c^{(1)}_{n-1}z^{n-1}+O(z^{n})\;$
where $c^{(1)}_0,\ldots,c^{(1)}_{n-1}$ are the numbers defined via 
\eqref{1.3}.} 

\medskip

Starting with a function $f_0:=f\in\cS$ of the form \eqref{1.1} and applying 
recursion \eqref{1.5} $n$ times one gets  a sequence of  Schur class functions 
satisfying
\begin{equation}
f_{j}(z)=\frac{zf_{j+1}(z)+c^{(j)}_0}{z\bar{c}^{(j)}_0f_{j+1}(z)+1}=
\frac{zf_{j+1}(z)+\gamma_j}{z\bar{\gamma}_jf_{j+1}(z)+1}\quad (j=0,\ldots,n)
\label{1.6}
\end{equation}
and such that $\; f_j(z)=c^{(j)}_0+c^{(j)}_1z+\ldots+c^{(j)}_{n-j}
z^{n-j}+O(z^{n-j+1})\;$ where $c^{(j)}_k$ are the numbers defined via 
\eqref{1.3}. Upon taking  the superposition of linear fractional 
transformations \eqref{1.6} one gets the linear fractional formula
\begin{equation}
f={\bf T}_\Theta[\cE]:=\frac{A\cE+B}{C\cE+D}
\label{1.7}
\end{equation}
which parametrizes all solutions to the ${\bf SP}_n$ where
the free parameter $\cE:=f_n$ runs through $\cS$ and the coefficient matrix
$\Theta={\scriptsize\left[\begin{array}{cc}A &B\\
C& D\end{array}\right]}$ is given by
\begin{equation}
\Theta(z)=W_0(z)W_1(z)\cdots W_n(z)\quad \mbox{where}\quad
W_j(z)=\left[\begin{array}{cc}z & \gamma_j \\ z\bar{\gamma}_j & 
1\end{array}\right]. 
\label{1.8}
\end{equation}
Motivated by engineering applications (where it is desirable for the 
solution $f$ of an interpolation problem to be rational and of small 
McMillan degree), the rational coefficient interpolation problem 
(as well as its multi-point analogs) was considered in  
\cite{aa} with an additional  constraint on the degree 
(complexity) of rational interpolants. In what follows, the polynomials
$N_f$ and $D_f$ will denote the numerator and the denominator from
the coprime representation $f=N_f/D_f$ of a rational function $f$.
By $\deg f=\max\{\deg N_f,\deg D_f\}$ we mean the McMillan degree of $f$.
The algebra of rational functions will be denoted by $\cR$ and we will let
$$
\cR_k:=\left\{f\in\cR: \; \deg f= k\right\}\quad\mbox{and}\quad
\cR_{\le k}:=\left\{f\in\cR: \; \deg f\le k\right\}.
$$
Being adapted to the single-point case, the problem formulated in \cite{aa} is:

\medskip
${\bf RP}_{n,k}$: {\em Given $c_0,\ldots,c_{n}\in\C$ and $k\ge 0$, find all
$f\in\cR_{\le k}$ of the form \eqref{1.1}.}

\medskip

The problem was solved in \cite{aa} and in \cite{abkw}
(for the matrix-valued case) as follows. 
\begin{theorem}
Let $q$ denote the rank of the Hankel matrix 
$H=\left[c_{i+j-1}\right]_{i,j\ge 1}$
constructed from the given numbers $c_j$ (the matrix $H$ is
$\frac{n-1}{2}\times \frac{n-1}{2}$ if $n$ is odd or $\frac{n-2}{2}\times
\frac{n}{2}$ if $n$ is even). Then
\begin{enumerate}
\item There is no $f\in\cR_k$ satisfying \eqref{1.1} for 
every $k<q$ or $q<k\le n-q$.
\item There exists  at most one function $f$ of complexity $k=q$ subject 
to \eqref{1.1}.
\item For every $k>n-q$, there are infinitely many solutions 
of the problem ${\bf RP}_{n,k}$ which are parametrized by the formula 
\begin{equation}
f={\bf T}_{\mathfrak A}[g]:=\frac{{\mathfrak
A}_{11}g+{\mathfrak A}_{12}}{{\mathfrak A}_{21}g+{\mathfrak A}_{22}},
\label{1.9}
\end{equation}
where the coefficients ${\mathfrak A}_{ij}$ are polynomials  explicitly  
constructed from the data set and such that 
$$
\deg\begin{bmatrix}{\mathfrak A}_{11} \\ {\mathfrak
A}_{21}\end{bmatrix}=q\quad\mbox{and}\quad
\deg\begin{bmatrix}{\mathfrak A}_{12} \\ {\mathfrak
A}_{22}\end{bmatrix}=n+1-q,
$$
and where the parameter $g=N_g/D_g\in\cR$ is such that 
$$
\deg N_g\le k-q,\quad \deg D_g\le k+q-n-1,\quad
{\mathfrak A}_{21}(0)N_g(0)+{\mathfrak A}_{22}(0)D_g(0)\neq 0.
$$
\end{enumerate}
\label{T:1.1}
\end{theorem}
We refer to \cite{abkw} for more details. In what follows, we use notation
$$
\cSR=\cS\cap \cR,\quad \cSR_{k}=\cS\cap \cR_k\quad \mbox{and}\quad\cSR_{\le
k}=\cS\cap \cR_{\le k}
$$
for the classes of functions in $\cR$, $\cR_k$ and $\cR_{\le k}$ respectively,
which are bounded by one in modulus on $\D$. Upon imposing both 
$H^\infty$-norm and complexity constraints (i.e., upon combining problems 
${\bf SP}_n$ and ${\bf RP}_{n,k}$) we arrive at the following interpolation 
problem.

\medskip
${\bf RSP}_{n,k}$: {\em Given an admissible data set $c_0,\ldots,c_{n}$ 
and $k\ge 0$, find all functions $f\in\cSR_{\le k}$ of the form \eqref{1.1}.}

\medskip

One may try to treat the latter problem using either formula \eqref{1.9} 
or \eqref{1.7}. In the first case, the complexity of $f$ is completely 
controlled by the complexity of the corresponding parameter $g$ and it 
suffices to pick up all parameters $g$ with $\deg g\le k-q$ leading via 
formula \eqref{1.9} to Schur-class functions $f$. However, this task is hard,
since formula  \eqref{1.9} does not control $\|{\bf T}_{\mathfrak 
A}[g]\|_\infty$ in terms of $\|g\|_\infty$. It may happen that a Schur 
class parameter $g$ 
produces $f\not\in\cS$ and on the other hand, a Schur class function  
$f\in\cSR_{\le k}$ may arise from a non-Schur class parameter $g$. 
Although Theorem \ref{T:1.1} guarantees that there are infinitely many 
functions $f\in\cR_{n+1-q}$ of the form \eqref{1.1}, it is not known
whether or not one of them is of the Schur class. The question about the 
minimal possible $k$ for which the problem ${\bf RSP}_{n,k}$ has a solution,
is still open. 

\medskip

It is not even clear from \eqref{1.9} that the problem ${\bf RSP}_{n,k}$ has 
solutions for $k$ large enough. On the other hand, the affirmative answer
for the latter question is readily seen from parametrization formula 
\eqref{1.7} which in contrast to \eqref{1.9}, perfectly controls the 
$H^\infty$-norm of $f$: all Schur-class rational solutions to the problem
${\bf SP}_n$ arise via formula \eqref{1.7} from some  Schur-class rational
parameter $\cE$. The complexities of interpolants are controlled here
to some extent. A straightforward induction argument deduces from
\eqref{1.8} that the coefficients $A$, $B$, $C$ and $D$ in \eqref{1.7} are 
polynomials of respective degrees $\deg A=n+1$, $\deg B\le n$, $\deg C\le 
n+1$,
$\deg D\le n$ and therefore,  
\begin{equation}
\deg {\bf T}_\Theta[\cE]\le n+1+\deg \cE.
\label{w}
\end{equation}
Letting $\cE$ in \eqref{1.7} to run through the class of constant functions 
(not exceeding one in modulus), one gets a family of solutions $f$ of the 
problem ${\bf RSP}_{n,n+1}$, but not all the solutions. It turns out that
zero cancellations may occur in \eqref{1.7} due to which some solutions 
to the ${\bf RSP}_{n,n+1}$ may arise from non-constant parameters. We also 
observe that the parameter $\cE\equiv 0$ leads via \eqref{1.7} to the function
${\bf T}_\Theta[0]=B/D\in\cSR_{\le n}$ which is therefore, a solution 
to the problem ${\bf RSP}_{n,n}$. The next example shows that this function 
might be {\em the only} solution to the  ${\bf RSP}_{n,n}$.
\begin{Ex}
{\rm Let $|c_0|<1$ and $c_j=0$ for $j=1,\ldots,n$. With this data, 
the problem ${\bf RSP}_{n,n}$ has only one solution $f\equiv c_0$. 
This follows from Theorem \ref{T:1.1} since in this case $q=0$.}
\label{E:1.2}
\end{Ex}
Otherwise (that is, if $c_j\neq 0$ at least for one $j\ge 1$ so that $q\ge 1$),
Theorem \ref{T:1.1} guarantees the existence of infinitely many functions
$f\in\cR_{\le n}$ of the form \eqref{1.1}, at least one of which (${\bf 
T}_\Theta[0]$)
belongs to $\cSR_{\le n}$. As was shown in \cite{blgm}-\cite{geor1}, the set 
of such 
functions is infinite and can be parametrized by
polynomials $\sigma$ with $\deg \sigma\le n$ and with all the roots outside $\D$.
More precisely, for every such $\sigma$, there exists a  unique (up to a common 
unimodular  constant factor) pair of polynomials $a(z)$ and $b(z)$, each of degree 
at most $n$  and such that
\begin{enumerate}
\item $|a(z)|^2-|b(z)|^2=|\sigma(z)|^2$ for $|z|=1$ and
\item the function $f=b/a$ (which belongs to $\cSR_{n}$ by part (1))
satisfies \eqref{1.1} and therefore, solves the ${\bf RSP}_{n,n}$.
\end{enumerate}
The objective of this note is to present an alternative parametrization of the 
solution set of the problem ${\bf RSP}_{n,k}$ (see Theorem \ref{T:1.3} below) which relies 
entirely on parametrization formula \eqref{1.7}. Some elementary analysis of the Schur 
algorithm will relate 
complexities of $\cE$ and $\deg {\bf T}_\Theta[\cE]$  more accurately than 
in \eqref{w}; this in turn, will allow us to describe all parameters 
$\cE\in\cSR$ leading via 
formula \eqref{1.7} to solutions $f$ of the problem ${\bf RSP}_{n,k}$ (these parameters will 
be called {\em admissible}). Explicit construction of these parameters 
is given below in terms of certain algorithm which seems to be quite efficient and simple 
from the computational point of view.  Here is the \underline{\bf Algorithm}:

\medskip

\noindent
{\bf Step 1:} {\em Given $c_0,\ldots,c_n$, compute the numbers
$\gamma_0,\gamma_1,\ldots,\gamma_n$ by
formula \eqref{1.4} using iteration \eqref{1.3}.}

\medskip

\noindent
{\bf Step 2:} {\em Using the numbers $\gamma_0,\ldots,\gamma_n$ compute
the polynomials 
\begin{equation}
A_{n}(z)=\sum_{j=0}^{n}a_jz^j \quad\mbox{and}\quad B_{n}(z)=\sum_{j=0}^{n}b_jz^j
\label{1.7a}
\end{equation}
from the system of recursions 
\begin{equation}
\left\{\begin{array}{l}A_0(z)\equiv \gamma_n,\quad B_0(z)\equiv 1,\\
\begin{array}{l}
A_{j+1}(z)=zA_{j}(z)+\gamma_{n-j-1}B_{j}(z),\\
B_{j+1}(z)=z\overline{\gamma}_{n-j-1}A_{j}(z)+B_{j}(z),\end{array}\quad
(j=0,\ldots,n-1).\end{array}\right.
\label{1.8a}
\end{equation}}
It is readily seen that $B_j(0)=1$ for $j=0,\ldots,n$. In particular,
$b_0=B_n(0)=1$.

\medskip

\noindent
{\bf Step 3:} {\em Using the coefficients $a_j$, $b_j$ from \eqref{1.7a} define 
the lower triangular Toeplitz matrices
\begin{equation}
{\bf A}=\left[\begin{array}{cccccc} a_{n}&0&\cdots&0\\ a_{n-1}
&a_{n}&\cdots&0\\ \vdots&\vdots&\ddots& 0\\ a_{1}&a_{2}&
\cdots&a_{n}\end{array}\right],\quad \widetilde{\bf B}=\left[\begin{array}{cccccc}
1&0&\cdots&0\\  \overline{b}_{1}
& 1&\cdots&0\\ \vdots&\vdots&\ddots &0\\
\overline{b}_{n-1}&\overline{b}_{n-2}&
\cdots&1\end{array}\right]
\label{1.9a}
\end{equation}
and compute the lower triangular Toeplitz matrix
\begin{equation}
{\bf R}=\M(r_{1},r_2\ldots,r_{n}):=\widetilde{\bf B}^{-1}{\bf A}.
\label{1.9b}
\end{equation}}
The three first steps are preliminary and can be carried out in finitely many steps.
The last step tells which parameters $\cE$ in \eqref{1.7} should be taken 
to get solutions to the problem ${\bf RSP}_{n,k}$. We first consider the 
case where $k=n$.

\medskip

\noindent
{\bf Step 4:} {\em For any $n$-tuple $\{\alpha_1,\ldots,\alpha_n\}$
of complex numbers, compute the function
\begin{equation}
\cE(z)=\frac{\beta_0+\beta_1z+\ldots +\beta_{n-1}z^{n-1}}
{\alpha_0+\alpha_1z+\ldots +\alpha_{n}z^{n}}
\label{1.10}
\end{equation}
where $\beta_0,\ldots,\beta_{n-1}$ are defined by
\begin{equation}
\left[\begin{array}{c}\beta_{n-1} \\ \vdots \\ \beta_0\end{array}\right]=
-{\bf R}\left[\begin{array}{c}\alpha_{n}\\ \vdots \\
\alpha_1\end{array}\right]
\label{1.11}
\end{equation}
where ${\bf R}$ is given in \eqref{1.9b} and $\alpha_0$ is such that 
$\cE\in\cS$}. 

\medskip

The main result of the paper is the following theorem; the proof will be 
given in Section 2.
\begin{Tm}
Let $\cE$ be constructed as in Step 4 and let $\Theta$ be as in \eqref{1.8}.
Then the function $f={\bf T}_\Theta[\cE]$ \eqref{1.7} solves the problem  
${\bf RSP}_{n,n}$ and conversely, all solutions of the ${\bf RSP}_{n,n}$ arise 
in this way. 
\label{T:1.3}
\end{Tm}
\begin{Rk}
\label{R:1.1}
{\rm The only relatively uncertain part in Step 4 is the choice of $\alpha_0$. 
However, it is readily seen that for any $\alpha_0$ satisfying
$|\alpha_0|\ge {\displaystyle\sum_{i=1}^{n}(|\alpha_i|+|\beta_{i-1}|)}$,
the function $\cE$ in \eqref{1.11} belongs to the Schur class which immediately
gives infinitely many solutions of the problem ${\bf RSP}_{n,n}$.
To be more precise, let us write \eqref{1.10} as 
$$
\cE(z)=\frac{P(z)}{\alpha_0+zQ(z)},
$$
where $P(z)=\beta_0+\beta_1z+\ldots +\beta_{n-1}z^{n-1}$
and $Q(z)=\alpha_1+\ldots +\alpha_{n}z^{n-1}$ and let $\D(c, r)$ denote  
the disk of radius $r$ centered at $c$. Then the set of all admissible 
$\alpha_0$'s
(for already chosen $\alpha_1,\ldots,\alpha_n$ and $\beta_0,\ldots,\beta_{n-1}$)
is the exterior (complement) of the domain $\Omega$ defined as
$$
\Omega=\bigcup_{|z|<1} \D(-zQ(z), |P(z)|).
$$}
\end{Rk}
\begin{Rk}
{\rm It follows from \eqref{1.10} that a parameter $\cE$ leading to a 
solution of the ${\bf RSP}_{n,n}$ has to satisfy $\cE(\infty)=0$. Thus, 
$\cE\equiv 0$ is 
the only admissible constant parameter for the problem ${\bf RSP}_{n,n}$ .  
Combining this fact with \eqref{w}, we conclude that every other constant function 
$\cE\in\cS$ leads via \eqref{1.7} to a solution of ${\bf RSP}_{n,n+1}$.}
\label{R:1.2}
\end{Rk}
As we have already seen, in contrast to the case $n=k$, 
the existence of infinitely many solutions of the problem ${\bf RSP}_{n,k}$ with 
$k>n$ is immediate. However, the description of all solutions is even somewhat 
more complicated. We get this description by an appropriate modification of 
Step 4 as follows.

\medskip

\noindent
{\bf Step 4$^\prime$:} {\em Let $k>n$ be fixed and let $\Theta$ and ${\bf R}$  
be as above. All solutions $f$ to the problem 
${\bf RSP}_{n,k}$ are obtained via formula \eqref{1.7} where the 
parameter $\cE$ is either any function from $\cSR_{\le k-n-1}$ or a 
function from $\cSR_{\le k}$ of the form 
\begin{equation}
\cE(z)=\frac{\beta_{n-k}+\beta_{n-k+1}z+\ldots +\beta_{n-1}z^{k-1}}
{\alpha_{n-k}+\alpha_{n-k+1}z+\ldots +\alpha_{n}z^{k}}
\label{1.10b}  
\end{equation}
where the coefficients  $\alpha_{n-k+1},\alpha_{n-k+2},\ldots,\alpha_n$ 
and  
$\beta_{n-k},\beta_{n-k+1},\ldots,\beta_{-1}$ are picked up arbitrarily,
after which the  coefficients $\beta_0,\ldots,\beta_{n-1}$ are 
defined as in \eqref{1.11} and where after all, the coefficient $\alpha_{n-k}$ 
is chosen so that the function $\cE$ of the form \eqref{1.10b} belongs to the 
Schur class $\cS$.}

\medskip

Justification of  Step 4$^\prime$ will be given in Section 2. In Section 3
we will present a version of Step 4 suitably modified for the case 
where $k<n$. There we will explain the reasons (by means of 
parametrization formula \eqref{1.7}) for which  the algorithm is not 
efficient for $k<n$.

\section{Proof of Theorem \ref{T:1.3}.}
\setcounter{equation}{0}

In this section we justify the algorithm presented in the previous section. Let 
\begin{equation}
\Theta_{k}(z):=W_{n-k}(z)\cdots W_n(z)
\label{2.0}
\end{equation}
where the factors $W_j$ are defined in \eqref{1.8}. 
Comparing \eqref{2.0} and \eqref{1.8} we see that $\Theta_n$
equals the coefficient matrix $\Theta$ of the transformation 
\eqref{1.7}. It is not hard to check by induction that $\Theta_{k}$ is of 
the form
\begin{equation}
\Theta_{k}(z)=\left[\begin{array}{cc}zB^\sharp_k(z) & A_k(z) \\ zA_k^\sharp(z) &
B_k(z)\end{array}\right]
\label{2.1}
\end{equation}
where  the polynomials $A_k$ and $B_k$ are constructed from system   
\eqref{1.8a} and where $A_k^\sharp$ and $B_k^\sharp$ are defined as follows:
\begin{equation}
A_k^\sharp(z)=z^k\overline{A_k(1/\bar{z})},\qquad
B_k^\sharp(z)=z^k\overline{B_k(1/\bar{z})}.
\label{2.1a}
\end{equation}
Let us take any $\cE=\frac{N_\cE}{D_\cE}\in\cSR$ and substitute it together 
with formula \eqref{2.1} for $\Theta_n=\Theta$ into \eqref{1.7}:
\begin{equation}
f(z)=\frac{zB^\sharp_n(z)N_\cE(z)+A_n(z)D_\cE(z)}
{zA^\sharp_n(z)N_\cE(z)+B_n(z)D_\cE(z)}.
\label{2.2}
\end{equation}
\begin{Rk}
{\rm The numerator and the denominator in \eqref{2.2} do not have common
zeros and thus,
\begin{equation}
N_{f}=zB^\sharp_nN_\cE+A_nD_\cE\quad\mbox{and}\quad
D_{f}=zA^\sharp_nN_\cE+B_nD_\cE.
\label{2.3}
\end{equation}}
\label{R:2.6}
\end{Rk}
{\bf Proof:} Taking determinants in \eqref{1.8}, \eqref{2.0} and 
\eqref{2.1} (with $k=n$) gives
\begin{eqnarray}
B_n(z)B_n^\sharp(z)-A_n(z)A_n^\sharp(z)&=&\frac{1}{z}\cdot\det\Theta_{n}(z)
\label{2.4}\\
&=&\frac{1}{z}\cdot\prod_{j=0}^n \det 
W_j(z)=z^{n}\cdot\prod_{j=0}^n(1-|\gamma_j|^2).\nonumber
\end{eqnarray}
Therefore, the only possible common zero for the numerator and 
the denominator
in \eqref{2.2} is $z=0$. But if this is the case, we then have 
$B_n(0)D_\cE(0)=D_\cE(0)=0$ which is impossible since the Schur function $\cE$ 
cannot have a pole at the origin.\qed

\medskip 

We shall now compare 
McMillan degrees of $f$ and $f_1$ in formula \eqref{1.5}.
\begin{La}
Let $f\in\cSR$ be of the form \eqref{1.5}. Then
$\deg f-1\le \deg f_1\le \deg f$. Moreover,
\begin{equation}
\deg f_1=\deg f  \; \; \Longleftrightarrow \; \; f_1(\infty)=0
\; \; \Longleftrightarrow \; \; f(\infty)\neq 1/\bar{c}_0
\label{2.5}
\end{equation}
and
\begin{equation}
\deg f_1=\deg f-1 \; \; \Longleftrightarrow \; \;  f_1(\infty)\neq 0 \; \;
\Longleftrightarrow \; \;  f(\infty)=1/\bar{c}_0.
\label{2.6}
\end{equation}
\label{L:2.2}
\end{La}
{\bf Proof:} Take $f_1$ in the form $f_1=N_{f_1}/D_{f_1}$ and
rewrite (\ref{1.5}) as
\begin{equation}
f(z)=\frac{zN_{f_1}(z)+c_0D_{f_1}(z)}{z\bar{c}_0N_{f_1}(z)+D_{f_1}(z)}=
\frac{F(z)}{G(z)}
\label{2.7}
\end{equation}  
from which we see that $\deg N_f\le  \deg f_1+1$,
$\deg D_f\le  \deg f_1+1$ and thus, $\deg f\le  \deg f_1+1$.
Now let us take $f$ in the form $f=\frac{N_f}{D_f}$ and solve equation
(\ref{1.5}) for $f_1$:
\begin{equation}
f_1(z)=\frac{(N_f(z)-c_0D_f(z))/z}{D_f(z)-\bar{c}_0N_f(z)}
\label{2.8}
\end{equation}
Since $c_0=f(0)=N_f(0)/D_f(0)$, it follows that
the numerator in \eqref{2.8} is a polynomial of degree not exceeding
$\deg f-1$. Therefore,
$\deg N_{f_1}\le  \deg f$, $\deg D_{f_1}\le  \deg \, f$ and thus,
$\deg f_1\le \deg f$. This completes the proof of the first
statement. 

\smallskip

Since there are only two possibilities for the value of 
$(\deg f-\deg f_1)$, statements \eqref{2.5} are equivalent to \eqref{2.6}.
We next observe that the  polynomials $F$ and $G$ in  \eqref{2.7} 
do not have common zeros (the proof is the same as in Lemma \ref{L:2.2})
and therefore we can conclude from \eqref{2.7} that
\begin{equation}
\deg f=\max\{\deg F, \, \deg G\}.
\label{2.9}
\end{equation}
Now we verify \eqref{2.5} (or (\ref{2.6})) separately for the following three
cases.

\medskip

\noindent
{\bf Case 1:} Let $\deg D_{f_1}>\deg N_{f_1}+1$. Then it follows from 
\eqref{2.7}
that $\deg f=\deg D_{f_1}=\deg
f_1$ and on the other hand, $f_1(\infty)=0$ and $f(\infty)=c_0\neq 1/\bar{c}_0$.

\medskip

\noindent
{\bf Case 2:} Let $\deg D_{f_1}<\deg N_{f_1}+1$. Then $\deg f=\deg 
N_{f_1}+1=\deg
f_1+1$ and on the other hand, $f_1(\infty)\neq 0$ and $f(\infty)=1/\bar{c}_0$.

\medskip

\noindent
{\bf Case 3:} Let $\deg D_{f_1}=\deg N_{f_1}+1$. Let $a_0$ and $b_0$ be the
leading coefficients of the polynomials $N_{f_1}$ and $D_{f_1}$ respectively.
Then the leading coefficients of $F$ and $G$ are $a_0+c_0b_0$ and 
$\bar{c}_0a_0+b_0=0$, respectively. Assuming that $\deg F< \deg 
N_{f_1}+1$ and $\deg G< \deg N_{f_1}+1$ we
have $a_0+c_0b_0=0$ and $\bar{c}_0a_0+b_0=0$ which gives
$a_0=b_0=0$ which is a contradiction. Therefore, 
$\max\{\deg F, \, \deg G\}=\deg N_{f_1}+1$ and by \eqref{2.9}, $\deg 
f=\deg N_{f_1}+1=\deg D_{f_1}=\deg f_1$.
Finally, since $\deg D_{f_1}=\deg N_{f_1}+1$, we have $f_1(\infty)=0$ and 
it follows from \eqref{2.7} that 
$f(\infty)={\displaystyle\frac{a_0+c_0b_0}{\bar{c}_0a_0+b_0}}$ which is not 
equal to $1/\bar{c}_0$, since $b_0\neq 0$ and $|c_0|\neq 1$.\qed

\medskip

Let us apply the backward Schur algorithm \eqref{1.6} to a function 
$\cE\in\cSR_k$ by letting
\begin{equation}
f_{n+1}=\cE\quad\mbox{and}\quad f_{j}(z)=
\frac{zf_{j+1}(z)+\gamma_j}{z\overline{\gamma}_jf_{j+1}(z)+1}
\quad\mbox{for} \quad j=0,\ldots,n.
\label{2.10}
\end{equation}
\begin{La}
If  $\deg f_{i}=\deg f_{i+1}+1$, then $\deg f_{j}=\deg f_{j+1}+1$ for every 
$j<i$.
If $f_i(\infty)=0$, then $f_j(\infty)=0$ and $\deg f_{j}=\deg f_{i}$ for every 
$j>i$.
\label{L:2.3}
\end{La}
{\bf Proof:} If $\deg f_{i}=\deg f_{i+1}+1$, then by virtue of \eqref{2.6}
(with $f$, $f_1$ and $c$ replaced respectively by  $f_i$, $f_{i+1}$ and 
$\gamma_i$) we have $f_i(\infty)=\frac{1}{\overline{\gamma}_i}\neq 0$. 
Then again by \eqref{2.6} (applied to the new triple $f_{i-1}$, $f_{i}$ 
and $\gamma_{i-1}$) we get $\deg f_{i-1}=\deg f_{i}+1$ and therefore,
$f_{i-1}(\infty)=\frac{1}{\overline{\gamma}_{i-1}}\neq 0$. The first statement
then follows by induction.

We now assume that $f_i(\infty)=0$. Since $f_i(\infty)\neq
\frac{1}{\overline{\gamma}_i}$, we conclude from \eqref{2.5} that
$f_{i+1}(\infty)=0$ and
$\deg f_{i+1}=\deg f_{i}$. The induction argument completes the proof of the 
second statement.\qed

\bigskip

\noindent
{\bf Proof of Theorem \ref{T:1.3}:} Let $f$ be a solution to the problem 
${\bf RSP}_{n,n}$, i.e., $f$ is a rational Schur-class function of 
degree at most $n$ satisfying equality \eqref{1.1}. Then $f$ is of the form 
\eqref{1.7} for some rational Schur-class function $\cE$ or equivalently, 
$f=f_0$ is obtained from $\cE=f_n$  via recursion \eqref{2.10}. Then we 
necessarily have
\begin{equation}
\deg \cE\le n\quad\mbox{and}\quad \cE(\infty)=0.
\label{2.10a}  
\end{equation}
Indeed, $\deg f\ge \deg \cE$ by Lemma \ref{L:2.2} and since 
$\deg f\le n$ by the assumption, the first relation in \eqref{2.10a} 
follows. If we assume that $\cE(\infty)\neq 0$, then we get by virtue of 
\eqref{2.6}, that $\deg f_{n-1}=\deg \cE+1$ and then we also have 
$\deg f=\deg \cE+n+1\ge n+1$ (by the first statement in Lemma \ref{L:2.3}) 
which contradicts the assumption. Thus, $\cE(\infty)=0$. Due to 
\eqref{2.10a} we can take $\cE$ in the form \eqref{1.10}, i.e., we can let
\begin{equation}
N_{\cE}(z)=\sum_{j=0}^{n-1}\beta_jz^j \quad\mbox{and}\quad 
D_{\cE}(z)=\sum_{j=0}^{n}\alpha_jz^j.
\label{2.11}
\end{equation}
It remains to show that the coefficients $\alpha_i$ and $\beta_i$ are 
related as in \eqref{1.11}. Observe, that the polynomials $A_n$ and $B_n$ 
constructed in \eqref{1.8a} are of degree at most $n$;
we take them in the form \eqref{1.7a} so that the reflected polynomials 
$A_{n}^\sharp$ and  $B_{n}^\sharp$ (see \eqref{2.1a}) are given by 
\begin{equation}
A_{n}^\sharp(z)=\sum_{j=0}^{n}\overline{a}_{n-j}z^j\quad\mbox{and}\quad
B_{n}^\sharp(z)=\sum_{j=0}^{n}\overline{b}_{n-j}z^j.
\label{2.12}
\end{equation}
Substituting \eqref{1.7a}, \eqref{2.11} and \eqref{2.12} into \eqref{2.3} 
we get
\begin{eqnarray}
N_{f}(z) 
&=&z^{n+1}\cdot\sum_{\ell=0}^{n-1}\left(\sum_{j=0}^{n-\ell-1}
(\overline{b}_{n-\ell-j-1}
\beta_{n-j-1}+a_{\ell+j+1}\alpha_{n-j})\right)z^\ell+P_1(z),\nonumber\\
D_{f}(z)&=&z^{n+1}\cdot\sum_{\ell=0}^{n-1}\left(\sum_{j=0}^{n-\ell-1}
(\overline{a}_{n-\ell-j-1}\beta_{n-j-1}+
b_{\ell+j+1}\alpha_{n-j})\right)z^\ell+P_2(z),\nonumber
\end{eqnarray}
where $P_{1}$ and $P_2$ are polynomials of degree at most $n$. The two 
latter formulas imply  that $\deg f\le n$ if and only if 
\begin{eqnarray}
\sum_{j=0}^{n-\ell-1}(\overline{b}_{n-\ell-j-1}\beta_{n-j-1}+
a_{\ell+j+1}\alpha_{n-j})=0\quad (\ell=0,\ldots,n-1),\label{2.13}\\
\sum_{j=0}^{n-\ell-1}(\overline{a}_{n-\ell-j-1}\beta_{n-j-1}+
b_{\ell+j+1}\alpha_{n-j})=0\quad (\ell=0,\ldots,n-1).
\label{2.14}
\end{eqnarray}
Making use of the Toeplitz matrices 
\begin{equation}
\begin{array}{ll}
{\bf A}=\M(a_1,a_2,\ldots,a_n),\quad & {\bf B}=\M(b_1,b_2,\ldots,b_n),\\
\widetilde{\bf A}=\M(\overline{a}_{n-1}, \overline{a}_{n-2},
\ldots,\overline{a}_0),\quad&
\widetilde{\bf B}=\M(\overline{b}_{n-1}, \overline{b}_{n-2}, 
\ldots,\overline{b}_{0}),\end{array}
\label{2.15}
\end{equation}
and of the vectors 
\begin{equation}
\bal=\left[\begin{array}{c}\alpha_n \\ \vdots \\ 
\alpha_1\end{array}\right]\quad\mbox{and}\quad
\bbe=\left[\begin{array}{c}\beta_{n-1} \\ \vdots \\ \beta_0\end{array}\right],
\label{2.16}
\end{equation}
one can write equations \eqref{2.13} and \eqref{2.14} in the matrix form as 
\begin{equation}
\widetilde{\bf B}\bbe+{\bf A}\bal=0\quad\mbox{and}\quad
\widetilde{\bf A}\bbe+{\bf B}\bal=0
\label{2.17}
\end{equation}
respectively. Since $b_0=B_n(0)=1$, the matrix $\widetilde{\bf B}$ is invertible.
Then we get from the first equation in \eqref{2.17} 
\begin{equation}
\bbe=-\widetilde{\bf B}^{-1}{\bf A}\bal=-{\bf R}\bal
\label{2.19}  
\end{equation}
which is the same as \eqref{1.11}. We thus showed that every solution 
$f$ to the problem ${\bf RSP}_{n,n}$ can be obtained via the Schur algorithm 
from a parameter $\cE\in\cS$ of the form \eqref{1.10}, \eqref{1.11}.

\smallskip

To show that  any such parameter is admissible, we have to verify that 
the vectors $\bal$ and  $\bbe$ related as in \eqref{2.19} satisfy both
equations in \eqref{2.17}. The first equation is clearly equivalent
to \eqref{2.19}. Substituting \eqref{2.19} into the second equation
and taking into account that all the matrices in \eqref{2.15} commute,
we get
\begin{equation}
\widetilde{\bf A}\bbe+{\bf B}\bal=-\widetilde{\bf A}\widetilde{\bf
B}^{-1}{\bf A}\bal+{\bf B}\bal=\widetilde{\bf B}^{-1}\left({\bf B}\widetilde{\bf
B}-{\bf A}\widetilde{\bf A}\right)\bal.
\label{2.19a}  
\end{equation}
We next substitute formulas \eqref{1.7a} and \eqref{2.12}
into \eqref{2.4} and examine the coefficients of $z^{2n-\ell}$ for 
$\ell=0,\ldots,n-1$ to get equalities
\begin{equation}
\sum_{j=0}^\ell(b_{n+j-\ell}\overline{b}_j-a_{n+j-\ell}\overline{a}_j)=0 \quad
(\ell=0,\ldots,n-1),
\label{2.18}
\end{equation}
which can be written in terms of matrices \eqref{2.15} as ${\bf 
B}\widetilde{\bf B}={\bf A}\widetilde{\bf A}$. We now conclude
from \eqref{2.19a} that the second equation in \eqref{2.17} is satisfied.
Thus, for every $\cE\in\cS$ of the form \eqref{1.10}, \eqref{1.11},
the coefficients $\alpha_i$, $\beta_i$ satisfy equalities \eqref{2.13},
\eqref{2.14} (i.e., equalities \eqref{2.17}), which in turn guarantees 
that the McMillan degree of the function $f$ obtained from $\cE$ via 
the Schur algorithm, does not exceed $n$. Since this $f$ belongs 
to $\cS$ and satisfies \eqref{1.1}, it solves the problem ${\bf RSP}_{n,n}$.\qed 

\bigskip

{\bf Justification of Step 4$^\prime$:} Let $k>n$ be a fixed integer. 
Every solution $f$ to the problem ${\bf RSP}_{n,k}$ is of the form \eqref{1.7} for 
some rational parameter $\cE\in\cSR$ with $\deg \cE \le k$. We have either 
$\cE(\infty)\neq 0$ or $\cE(\infty)=0$. In the first case, $\deg f=\deg\cE+n+1$ 
(by Lemmas \ref{L:2.2} and \ref{L:2.3}) and therefore, $\deg \cE \le k-n-1$. On the 
other hand, for every $\cE\in\cSR_{\le k-n-1}$, it follows from \eqref{w} that 
$\deg {\bf T}_\Theta[\cE]\le k$. In the second case, we can take $\cE$ in the 
form \eqref{1.10b}, that is to let  
$$
N_{\cE}(z)=\sum_{j=0}^{k-1}\beta_{n-k+j}z^j \quad\mbox{and}\quad
D_{\cE}(z)=\sum_{j=0}^{k}\alpha_{n-k+j}z^j.
$$
Substituting the latter formulas along with \eqref{1.7a} and \eqref{2.12} into 
\eqref{2.3} we get the formulas for $N_f$ and $D_f$ as in the proof of Theorem 
\ref{T:1.3} but with the factor $z^{k+1}$ (rather than $z^{n+1}$)
on the left and with 
polynomials $P_1$ and $P_2$ of degree at most $k$. Then we conclude that $\deg 
f\le k$ if and only if conditions \eqref{2.17} hold which is equivalent to 
\eqref{2.19}.

\section{Concluding remarks} 
\setcounter{equation}{0}

In conclusion we present a version of the main algorithm for the case where $k<n$.
The three first steps are the same as before; the last step describing all admissible 
parameters in parametrization formula \eqref{1.7} is the following.

\medskip

\noindent
{\bf Step 4$^{\prime\prime}$:} {\em Let $k<n$ be fixed and let $\Theta$ and 
${\bf R}=\M(r_{1}, r_{2},\ldots,r_{n})$ be as above. All 
solutions $f$ 
to the problem ${\bf RSP}_{n,k}$ are obtained via formula \eqref{1.7} where the
parameter $\cE$ is a Schur-class function of the form
\begin{equation}
\cE(z)=\frac{N_\cE(z)}{D_\cE(z)}=\frac{\beta_0+\beta_1z+\ldots 
+\beta_{k-1}z^{k-1}}
{\alpha_0+\alpha_1z+\ldots +\alpha_{k}z^{k}}.
\label{4.1}   
\end{equation}
where the coefficients  $\alpha_0,\ldots,\alpha_k$
and $\beta_0,\ldots,\beta_{k-1}$ satisfy the system
\begin{equation}
\left[\begin{array}{c}\beta_{k-1} \\ \vdots \\ \beta_1 \\ 
\beta_0\end{array}\right]=
\left[\begin{array}{cccccc}{r}_{n}&0&\cdots&0\\ r_{n-1}
& r_{n}&\cdots&0\\ \vdots&\vdots&\ddots &0\\r_{n-k+1}&r_{n-k}&
\cdots&r_{n}\end{array}\right]\left[\begin{array}{c}\alpha_{k} \\ \vdots \\
\alpha_2\\ \alpha_1\end{array}\right],\quad
\label{4.5}
\end{equation}
\begin{equation}
\left[\begin{array}{cccc}r_1 & r_2 & \ldots & r_{k+1} \\
r_2 & r_3 & \ldots & r_{k+2}\\ \vdots & \vdots && \vdots \\
r_{n-k} & r_{n-k+1} & \ldots &
r_{n}\end{array}\right]\left[\begin{array}{c}\alpha_{k} \\
\alpha _{k-1} \\ \vdots \\ \alpha_0\end{array}\right]=0,
\label{4.6}
\end{equation}
(the matrix in \eqref{4.6} is of Hankel structure).}

\medskip

{\bf Proof:} As in the proof of Theorem \ref{T:1.3} we first observe that 
every solution $f$ of the problem ${\bf RSP}_{n,k}$ is of the form 
\eqref{1.7} for some $\cE\in\cSR_{\le k}$ subject to $\cE(\infty)=0$.
Therefore, $\cE$ can be taken in the form \eqref{4.1}.
Substituting \eqref{1.7a}, \eqref{2.12} and \eqref{4.1} into \eqref{2.3} 
we now get $N_f$ and $D_f$ the polynomials of degree at most $n+k$.
Then equating the coefficients of $z^{k+\ell}$ of these polynomials to 
zero for $j=\ell,\ldots,n-1$, we get necessary and sufficient conditions 
(similar to \eqref{2.13} and \eqref{2.14}) for 
$\deg f=\max\{\deg N_f,\deg D_f\}$ not to exceed $k$. These conditions are
\begin{eqnarray}
\sum_{j=0}^{{\rm min}\{n-\ell-1,k-1\}}\overline{b}_{n-\ell-j-1}\beta_{k-j-1}+
\sum_{j=0}^{{\rm min}\{n-\ell-1,k\}}a_{\ell+j+1}\alpha_{k-j}=0,\nonumber\\
\sum_{j=0}^{{\rm 
min}\{n-\ell-1,k-1\}}\overline{a}_{n-\ell-j-1}\beta_{k-j-1}+
\sum_{j=0}^{{\rm min}\{n-\ell-1,k\}}b_{\ell+j+1}\alpha_{k-j}=0\nonumber  
\end{eqnarray}
($\ell=0,\ldots,n$) and it is not hard to see that they can be written in 
the matrix form \eqref{2.17} as 
\begin{equation}
\widetilde{\bf B}\bbe+{\bf A}\bal=0\quad\mbox{and}\quad
\widetilde{\bf A}\bbe+{\bf B}\bal=0
\label{4.2}
\end{equation}
respectively where the matrices ${\bf A}$, ${\bf B}$, $\widetilde{\bf A}$ 
and $\widetilde{\bf B}$ are the same as in \eqref{2.15} and where now  
\begin{equation}
\bal=\left[\begin{array}{c}\alpha_{k} \\ \vdots \\ \alpha_1 \\
\alpha_0 \\ 0 \\ \vdots \\ 0\end{array}\right]\quad\mbox{and}\quad
\bbe=\left[\begin{array}{c}\beta_{k-1} \\ \vdots \\ \beta_0 \\ 0\\ 0 \\ \vdots
\\ 0\end{array}\right].
\label{4.3}
\end{equation}
Since ${\bf B}\widetilde{\bf B}={\bf A}\widetilde{\bf A}$, it follows 
as in the proof of Theorem \ref{T:1.3}, that $\bal$ and $\bbe$ solve the 
system \eqref{4.2} if and only if  they are related as in \eqref{2.19}. 
Substituting \eqref{4.3} into \eqref{2.19} and comparing the $k$ top 
entries in the obtained equality, we get \eqref{4.5}; comparison 
of the $n-k$ bottom entries gives \eqref{4.6}. \qed
\begin{Rk}
{\rm Although Step 4$^{\prime\prime}$ looks very similar to Step 4 in 
Section 1, in fact it is much less efficient. Let us demonstrate this by the 
case where $k=n-1$. Then condition \eqref{4.6} takes the form 
\begin{equation}
r_1\alpha_{n-1}+r_2\alpha_{n-2}+\ldots+r_{n-1}\alpha_1+r_{n}\alpha_0=0
\label{4.7}
\end{equation}
and if $r_{n}\neq 0$, then $a_0$ is uniquely determined by 
$\alpha_1,\ldots,\alpha_{n-1}$. The problem is to describe all the
tuples $\{\alpha_1,\ldots,\alpha_{n-1}\}$ (which now are the only free 
parameters) for which the function 
\begin{equation}
\cE(z)=\frac{\beta_0+\beta_1z+\ldots
+\beta_{n-2}z^{n-2}}
{\alpha_0+\alpha_1z+\ldots +\alpha_{n-1}z^{n-1}}
\label{4.8}   
\end{equation}
with the coefficients $\alpha_0$, $\beta_0,\ldots,\beta_{n-2}$ determined by 
formulas \eqref{4.7} and \eqref{4.5} (with $k=n-1$), belongs to the Schur class.
The problem is hard; at the moment we even do not know necessary and sufficient 
conditions  for the existence of at least one such tuple.} 
\label{R:3.2}
\end{Rk}

\bibliographystyle{amsplain}
\providecommand{\bysame}{\leavevmode\hbox to3em{\hrulefill}\thinspace}

\end{document}